\documentclass[11pt]{article}
\usepackage{mathrsfs,amsmath,amssymb,epsfig,subfigure,color}

\renewcommand{\arraystretch}{1.2}
\newcommand{\IR}{{\mathbb{R}}}

\newcommand{\BE}{\begin{equation}}
\newcommand{\EE}{\end{equation}}

\newtheorem{theorem}{Theorem}[section]
\newtheorem{lemma}{Lemma}[section]
\newtheorem{example}{Example}[section]

\newtheorem{remark}{Remark} [section]

\def\dfrac{\displaystyle\frac}
\def\dsum{\displaystyle\sum}
\def\dint{\displaystyle\int}

\begin{document}
\title{\textbf{High order difference schemes for a time fractional differential equation with Neumann boundary conditions}}
\author{Seakweng Vong\thanks{Email: swvong@umac.mo. Department of Mathematics, University of Macau, Av. Padre Tom\'{a}s Pereira Taipa, Macau, China.}
\and Zhibo Wang\thanks{Corresponding author. Email: zhibowangok@gmail.com. Department of Mathematics, University of Macau, Av. Padre Tom\'{a}s Pereira Taipa, Macau, China.}}
\date{}
 \maketitle
\begin{abstract}
 Based on our recent results, in this paper,
 a compact finite difference scheme is derived for a time fractional differential equation
 subject to the Neumann boundary conditions.
 The proposed scheme is second order accurate in time and fourth order accurate in space.
 In addition, a high order alternating direction implicit (ADI) scheme is also constructed for the two-dimensional case. Stability and convergence of the schemes are analyzed using their matrix forms.
\end{abstract}

\section{Introduction}
 Fractional differential equations have grown to be the focus of many studies due to their various applications.
 Readers can refer to the books \cite{Podlubny,Kilbas} for background of these equations.
 One of the key features of the fractional derivatives is the nonlocal dependence which causes difficulties when numerical schemes
 are designed for solving fractional differential equations. However, with
 the efforts of numerous researchers, great progress along this direction
 has been made in recent years. Interested readers can refer to
 \cite{Lubich}--\cite{Rensci} for a brief review. We remark here that
 the list does not mean to be complete but we try to include those that are more
 related to the present study.

 In this article, we consider high order finite difference schemes for the following time fractional differential equation in a region $\Omega$:
 \begin{equation}\label{original}
 _{0}^CD_t^{\gamma}u({\bf x},t)=\kappa_1\Delta u({\bf x},t)
 -\kappa_2u({\bf x},t)+g({\bf x},t),
 \quad{\bf x}\in\Omega, \quad 0<t\leq T, \quad 1<\gamma<2,\\
 \end{equation}
 subject to the initial conditions:
 $$u({\bf x},0)=\psi({\bf x}), \quad \dfrac{\partial u({\bf x},0)}{\partial t}=\phi({\bf x}), \quad {\bf x}\in\bar{\Omega}=\Omega\cup\partial\Omega,$$
 and the zero flux boundary condition:
 $$\frac{\partial u({\bf x},t)}{\partial {\bf n}}=0,
 \quad {\bf x}\in \partial\Omega,\quad 0<t\leq T,$$
 where $\partial \Omega$ is the boundary of $\Omega$, $\frac{\partial}{\partial {\bf n}}$ is the
 differentiation in the normal direction
 and $\kappa_1$, $\kappa_2$ are some positive constants. We further suppose that $\frac{\partial \psi({\bf x})}{\partial {\bf n}}=0,$ for ${\bf x}\in \partial\Omega$.
 We have used $_{0}^CD_t^{\gamma}u$ to denote the Caputo
 fractional derivative of $u$ with respect to the time variable $t$, which is $$_{0}^CD_t^{\gamma}u({\bf x},t)=\frac1{\Gamma(2-\gamma)}\int_0^t
 \frac{\partial^2u({\bf x},s)}{\partial s^2}(t-s)^{1-\gamma}ds,$$
 with $\Gamma(\cdot)$ being the gamma function. Theoretical results such as existence and uniqueness of solutions
 to fractional differential equations can be found in \cite{Podlubny,Kilbas}.
 In recent years, there are growing interests on the study of numerical solutions
 for time fractional differential equations subject to the Neumann boundary condition \cite{Yuste1}--\cite{Rensci}.

 We note that equation (\ref{original}) can be written equivalently as \cite{Huang}:
 $$\dfrac{\partial u({\bf x},t)}{\partial t}=\phi({\bf x})+
 \dfrac1{\Gamma(\alpha)}\dint_0^t(t-s)^{\alpha-1}
 \bigg[\kappa_1\dfrac{\partial^2u({\bf x},s)}{\partial x^2}-\kappa_2u({\bf x},s)\bigg]ds+f({\bf x},t),\quad{\bf x}\in\Omega, \quad 0<t\leq T,$$
 where $0<\alpha=\gamma-1<1,~f({\bf x},t)={_{0}I^\alpha_tg({\bf x},t)}$, and $_{0}I^\alpha_t$ is the Riemann-Liouville fractional integral operator of order $\alpha$, defined as
 $$_{0}I^\alpha_tg({\bf x},t)
 =\dfrac1{\Gamma(\alpha)}\dint_0^t(t-s)^{\alpha-1}g({\bf x},s)ds.$$

 By applying the weighted and shifted Gr\"{u}nwald difference
 (see \cite{Tian,Zhou, Wang}) to the
 Riemann-Liouville fractional integral, we establish compact schemes with
 second order temporal accuracy and fourth order spatial accuracy.
 Our analysis is based on the matrix form of the schemes and it turns out to give
 intuitive ideas of some norms and inner products defined in previous related works.

 This paper is organized as follows. We first consider the one-dimensional problem in Section 2 and 3,
 where we propose a high order scheme and study its convergence respectively.
 In Section 4, a high order alternating direction implicit scheme is proposed for the two-dimensional problem.
 Numerical examples are given in the last section.

 \section{The proposed compact difference scheme}
 In this section, we develop a high order scheme for the following one-dimensional problem:
 \begin{align}\label{main1}
 &_{0}^CD_t^{\gamma}u(x,t)
 =\kappa_1\dfrac{\partial^2u(x,t)}{\partial x^2}-\kappa_2u(x,t)+g(x,t),
  \quad 0\leq x\leq L, \quad 0<t\leq T, \quad 1<\gamma<2,\\\label{main2}
 &u(x,0)=\psi(x), \quad \dfrac{\partial u(x,0)}{\partial t}=\phi(x), \quad 0\leq x\leq L,\\\label{main3}
 &\frac{\partial u(0,t)}{\partial x}=0, \quad \frac{\partial u(L,t)}{\partial x}=0, \quad 0<t\leq T.
 \end{align}
 We assume that $\psi\equiv0$ in (\ref{main2}) without loss of generality
 since we can solve the equation for $v(x,t)=u(x,t)-\psi(x)$ in general.

 An equivalent form of (\ref{main1}) can be read as:
 \begin{equation}\label{changed-main1}
 \dfrac{\partial u(x,t)}{\partial t}=\phi(x)+
 \dfrac1{\Gamma(\alpha)}\dint_0^t(t-s)^{\alpha-1}
 \bigg[\kappa_1\dfrac{\partial^2u(x,s)}{\partial x^2}-\kappa_2u(x,s)\bigg]ds+f(x,t),
 \end{equation}
 where $0\leq x\leq L,~0<t\leq T,~0<\alpha=\gamma-1<1,~f(x,t)={_{0}I^\alpha_tg(x,t)}$.

 For given integers $M$ and $N$, we discretize the equation using spatial step size $h=\frac LM$ and temporal step size $\tau=\frac TN$ respectively.
 For $i=0,1,\ldots,M$ and $k=0,1,\ldots,N$, denote $x_i=ih,~t_k=k\tau,~u^k=(u_0^k,u_1^k,\ldots,u_M^k)^T$.
 To study the grid function $u=\{u_i^k|0\leq i\leq M,~ 0\leq k\leq N\}$ that approximates the solution, the following notations are needed:
 $$\delta_xu_{i-\frac12}^k=\dfrac1h(u_i^k-u_{i-1}^k),~~1\leq i\leq M,$$
 $$\delta_x^2u_{i}^k=\left\{
 \begin{array}{ll}
 \frac2h\delta_xu_{\frac12}^k,&i=0,\\
 \frac1h(\delta_xu_{i+\frac12}^k-\delta_xu_{i-\frac12}^k),&1\leq i\leq M-1,\\
 -\frac2h\delta_xu_{M-\frac12}^k,&i=M.
 \end{array}\right.$$
 $${\cal H}u_i=\left\{
 \begin{array}{ll}
 \frac16(5u_0+u_1), &i=0,\\
 \frac1{12}(u_{i-1}+10u_i+u_{i+1}), &1\leq i\leq M-1,\\
 \frac16(u_{M-1}+5u_M), &i=M
 \end{array}\right.$$
 $$\langle u,v\rangle=h\dsum_{i=0}^{M}u_iv_i,~~\|u\|^2=\langle u,u\rangle,~~\|u\|_\infty=\max_{0\leq i\leq M}|u_i|.$$
 Discretization of $\frac{\partial^2 u}{\partial x^2}$ is based on the following lemma:
 \begin{lemma}\label{compact}{\rm(\cite{Rensci})}
 Denote $\zeta(s)=(1-s)^3[5-3(1-s)^2]$.

 {\rm(I)} If $f(x)\in {\cal C}^6[x_0,x_1]$, then we have
 $$\begin{array}{rl}
 &\quad\Big[\dfrac56f''(x_0)+\dfrac16f''(x_1)\Big]-\dfrac2h\Big[\dfrac{f(x_1)-f(x_0)}h-f'(x_0)\Big]\\[5pt]
 &=-\dfrac h6f'''(x_0)+\dfrac{h^3}{90}f^{(5)}(x_0)+
 \dfrac{h^4}{180}\dint_0^1f^{(6)}(x_0+sh)\zeta(s)ds.
 \end{array}$$

 {\rm(II)} If $f(x)\in {\cal C}^6[x_{M-1},x_M]$, then we get
 $$\begin{array}{rl}
 &\quad\Big[\dfrac16f''(x_{M-1})+\dfrac56f''(x_M)\Big]-\dfrac2h\Big[f'(x_M)-\dfrac{f(x_M)-f(x_{M-1})}h\Big]\\[5pt]
 &=\dfrac h6f'''(x_M)-\dfrac{h^3}{90}f^{(5)}(x_M)+
 \dfrac{h^4}{180}\dint_0^1f^{(6)}(x_M-sh)\zeta(s)ds.
 \end{array}$$

 {\rm(III)} If $f(x)\in {\cal C}^6[x_{i-1},x_{i+1}],~1\leq i\leq M-1$, then it holds that
 $$\begin{array}{rl}
 &\quad\dfrac1{12}[f''(x_{i-1})+10f''(x_i)+f''(x_{i+1})]
 -\dfrac1{h^2}[f(x_{i-1})-2f(x_i)+f(x_{i+1})]\\[5pt]
 &=\dfrac{h^4}{360}\dint_0^1[f^{(6)}(x_i-sh)+f^{(6)}(x_i+sh)]\zeta(s)ds.
 \end{array}$$
 \end{lemma}

 As mentioned before, our scheme for equation (\ref{main1}) is derived using the equivalent form (\ref{changed-main1}). We need to introduce the following shifted Gr\"{u}nwald difference to the Riemann-Liouville fractional integral:
 $$\mathcal{A}_{\tau,r}^\alpha f(t)=\tau^\alpha\dsum_{k=0}^{\infty}\omega_kf(t-(k-r)\tau),$$
 where $\omega_k=(-1)^k\binom{-\alpha}{k}$.
 By the idea of \cite{Tian}, the following second order approximation for Riemann-Liouville fractional integrals is derived very recently in \cite{Wang}.
 \begin{lemma}\label{main-operator}
 Let $f(t),~_{-\infty}I_t^{2-\alpha}f$ and $(i\omega)^{2-\alpha}\mathscr{F}[f](\omega)$ belong to $L^1(\mathbb{R})$. Define the weighted and shifted difference operator by
 $$\mathcal{I}_{\tau,p,q}^\alpha f(t)
 =\dfrac{2q+\alpha}{2(q-p)}\mathcal{A}_{\tau,p}^\alpha f(t)
 +\dfrac{2p+\alpha}{2(p-q)}\mathcal{A}_{\tau,q}^\alpha f(t),$$
 then we have
 $$\mathcal{I}_{\tau,p,q}^\alpha f(t)={_{-\infty}I_t^\alpha f(t)}+O(\tau^2)$$
 for $t\in \mathbb{R}$, where $p$ and $q$ are integers and $p\neq q$.
 \end{lemma}

 With $(p,q)=(0,-1)$, which yields $\frac{2q+\alpha}{2(q-p)}=1-\frac\alpha2,~\frac{2p+\alpha}{2(p-q)}=\frac\alpha2$
 in Lemma \ref{main-operator}, we get that
 $$\begin{array}{rl}
 {_{0}I^\alpha_tu(x_i,t_{n+1})}&=\tau^\alpha\bigg[\Big(1-\dfrac\alpha2\Big)\dsum_{k=0}^{n+1}\omega_ku^{n+1-k}_i+
 \dfrac\alpha2\dsum_{k=0}^{n}\omega_ku^{n-k}_i\bigg]+O(\tau^2)\\
 &=\tau^\alpha\dsum_{k=0}^{n+1}\lambda_ku^{n+1-k}_i+O(\tau^2),
 \end{array}$$
 and
 $$\begin{array}{rl}
 {_{0}I^\alpha_tu_{xx}(x_i,t_{n+1})}&=\tau^\alpha\bigg[\Big(1-\dfrac\alpha2\Big)\dsum_{k=0}^{n+1}\omega_k\delta^2_xu^{n+1-k}_i+
 \dfrac\alpha2\dsum_{k=0}^{n}\omega_k\delta^2_xu^{n-k}_i\bigg]+O(\tau^2+h^2)\\
 &=\tau^\alpha\dsum_{k=0}^{n+1}\lambda_k\delta^2_xu^{n+1-k}_i+O(\tau^2+h^2),
 \end{array}$$
  where
 \begin{equation}\label{sequence-lambda}
 \lambda_0=(1-\dfrac\alpha2)\omega_0,
 ~\lambda_k=(1-\dfrac\alpha2)\omega_{k}+\dfrac\alpha2\omega_{k-1},~k\geq1.
 \end{equation}
 Therefore, a weighted Crank-Nicolson scheme for equation (\ref{changed-main1}) can be given by
 $$\dfrac{u^{n+1}_i-u^n_i}{\tau}=\phi_i+\dfrac{\tau^{\alpha}}2
 \bigg[\dsum_{k=0}^{n+1}\lambda_k(\kappa_1\delta^2_xu^{n+1-k}_i-\kappa_2u^{n+1-k}_i)
 +\dsum_{k=0}^{n}\lambda_k(\kappa_1\delta^2_xu^{n-k}_i-\kappa_2u^{n-k}_i)\bigg]
 +\dfrac12(f^n_i+f^{n+1}_i).$$
 To derive a higher order scheme, we follow the idea in \cite{Rensci}.
 Beginning with $i=0$,  one has
 \begin{equation}\label{compact-0}
 \begin{array}{ll}
 \quad{\cal H}(u^{n+1}_0-u^n_0)\\
 =\dfrac{\tau^{\alpha+1}}2
 \dsum_{k=0}^{n+1}\lambda_k\bigg(\kappa_1\Big[\frac2h\delta_xu_{\frac12}^{n+1-k}
 -\frac2h\frac{\partial u(0,t_{n+1-k})}{\partial x}
 -\frac h6\frac{\partial^3u(0,t_{n+1-k})}{\partial x^3}
 +\frac {h^3}{90}\frac{\partial^5u(0,t_{n+1-k})}{\partial x^5}\Big]
 -\kappa_2{\cal H}u^{n+1-k}_0\bigg)\\
 \quad+\dfrac{\tau^{\alpha+1}}2
 \dsum_{k=0}^{n}\lambda_k\bigg(\kappa_1\Big[\frac2h\delta_xu_{\frac12}^{n-k}
 -\frac2h\frac{\partial u(0,t_{n-k})}{\partial x}
 -\frac h6\frac{\partial^3u(0,t_{n-k})}{\partial x^3}
 +\frac {h^3}{90}\frac{\partial^5u(0,t_{n-k})}{\partial x^5}\Big]
 -\kappa_2{\cal H}u^{n-k}_0\bigg)\\
 \quad+\tau{\cal H}\phi_0+\dfrac{\tau}2{\cal H}(f^n_0+f^{n+1}_0)
 +\tau R_0^{n+1},
 \end{array}
 \end{equation}
 where $R_0^{n+1}=O(\tau^2+h^4)$.

 \smallskip
  We can now differentiate equation (\ref{main1}) with respect to $x$ to give
 $$_{0}^CD_t^{\gamma}\frac{\partial u(x,t)}{\partial x}
 =\kappa_1\dfrac{\partial^3u(x,t)}{\partial x^3}-\kappa_2\frac{\partial u(x,t)}{\partial x}+g_x(x,t).$$
 Letting $x\rightarrow0^+$ and noticing the boundary condition (\ref{main3}), we have
 \begin{equation}\label{skill1}
 \kappa_1\dfrac{\partial^3u(0,t)}{\partial x^3}=-g_x(0,t).
 \end{equation}
 With the Caputo fractional derivative operator $_{0}^CD_t^{\gamma}$ acting on (\ref{skill1}),  it follows that
 \begin{equation}\label{skill2}
 {_{0}^CD_t^{\gamma}}\dfrac{\partial^3u(0,t)}{\partial x^3}
 =-\frac1{\kappa_1}{_{0}^CD_t^{\gamma}}g_x(0,t).
 \end{equation}
 Meanwhile, differentiating equation (\ref{main1}) three times with respect to $x$ yields
 \begin{equation}\label{skill3}
 _{0}^CD_t^{\gamma}\frac{\partial^3u(x,t)}{\partial x^3}
 =\kappa_1\dfrac{\partial^5u(x,t)}{\partial x^5}
 -\kappa_2\frac{\partial^3u(x,t)}{\partial x^3}+g_{xxx}(x,t).
 \end{equation}
 Once again, let $x\rightarrow0^+$ in (\ref{skill3}). We can then substitute (\ref{skill1}) and (\ref{skill2}) to (\ref{skill3}) to achieve
 \begin{equation}\label{skill4}
 \kappa_1\dfrac{\partial^5u(0,t)}{\partial x^5}=-g_{xxx}(0,t)-\frac{\kappa_2}{\kappa_1}g_x(0,t)-\frac1{\kappa_1}{_{0}^CD_t^{\gamma}g_x(0,t)}.
 \end{equation}
 Inserting (\ref{skill1}),~(\ref{skill4}) into (\ref{compact-0})
 and noticing the boundary condition (\ref{main3}), the compact scheme for $i=0$ can be given,
 by omitting small terms, as:
 \begin{equation}\label{compact-scheme1}
 \begin{array}{ll}
 \quad{\cal H}(u^{n+1}_0-u^n_0)\\
 =\dfrac{\tau^{\alpha+1}}2
 \bigg[\dsum_{k=0}^{n+1}\lambda_k(\kappa_1\delta^2_xu^{n+1-k}_0
 -\kappa_2{\cal H}u^{n+1-k}_0)
 +\dsum_{k=0}^{n}\lambda_k(\kappa_1\delta^2_xu^{n-k}_0
 -\kappa_2{\cal H}u^{n-k}_0)\bigg]\\
 \quad+\dfrac{\tau^{\alpha+1}}2\dsum_{k=0}^{n+1}
 \lambda_k\bigg(\frac h6(g_x)_0^{n+1-k}-\frac{h^3}{90}
 \Big[(g_{xxx})_0^{n+1-k}+\frac{\kappa_2}{\kappa_1}(g_x)_0^{n+1-k}+\frac1{\kappa_1}({_{0}^CD_t^{\alpha+1}g_x})_0^{n+1-k}\Big]\bigg)\\
 \quad+\dfrac{\tau^{\alpha+1}}2\dsum_{k=0}^{n}
 \lambda_k\bigg(\frac h6(g_x)_0^{n-k}-\frac{h^3}{90}
 \Big[(g_{xxx})_0^{n-k}+\frac{\kappa_2}{\kappa_1}(g_x)_0^{n-k}+\frac1{\kappa_1}({_{0}^CD_t^{\alpha+1}g_x})_0^{n-k}\Big]\bigg)\\
 \quad+\tau{\cal H}\phi_0+\dfrac\tau2{\cal H}(f^n_0+f^{n+1}_0),\quad 0\leq n\leq N-1.
 \end{array}
 \end{equation}
 The scheme at the other end can be similarly derived as
 \begin{equation}\label{compact-scheme3}
 \begin{array}{ll}
 \quad{\cal H}(u^{n+1}_M-u^n_M)\\
 =\dfrac{\tau^{\alpha+1}}2
 \bigg[\dsum_{k=0}^{n+1}\lambda_k(\kappa_1\delta^2_xu^{n+1-k}_M
 -\kappa_2{\cal H}u^{n+1-k}_M)
 +\dsum_{k=0}^{n}\lambda_k(\kappa_1\delta^2_xu^{n-k}_M
 -\kappa_2{\cal H}u^{n-k}_M)\bigg]\\
 \quad-\dfrac{\tau^{\alpha+1}}2\dsum_{k=0}^{n+1}
 \lambda_k\bigg(\frac h6(g_x)_M^{n+1-k}-\frac{h^3}{90}
 \Big[(g_{xxx})_M^{n+1-k}+\frac{\kappa_2}{\kappa_1}(g_x)_M^{n+1-k}+\frac1{\kappa_1}({_{0}^CD_t^{\alpha+1}g_x})_M^{n+1-k}\Big]\bigg)\\
 \quad-\dfrac{\tau^{\alpha+1}}2\dsum_{k=0}^{n}
 \lambda_k\bigg(\frac h6(g_x)_M^{n-k}-\frac{h^3}{90}
 \Big[(g_{xxx})_M^{n-k}+\frac{\kappa_2}{\kappa_1}(g_x)_M^{n-k}+\frac1{\kappa_1}({_{0}^CD_t^{\alpha+1}g_x})_M^{n-k}\Big]\bigg)\\
 \quad+\tau{\cal H}\phi_M+\dfrac\tau2{\cal H}(f^n_M+f^{n+1}_M),\quad 0\leq n\leq N-1.
 \end{array}
 \end{equation}
 One can readily see that, at the internal grid, the scheme can be written as
 \begin{equation}\label{compact-scheme2}
 \begin{array}{ll}
 \quad{\cal H}(u^{n+1}_i-u^n_i)\\
 =\dfrac{\tau^{\alpha+1}}2
 \bigg[\dsum_{k=0}^{n+1}\lambda_k(\kappa_1\delta^2_xu^{n+1-k}_i
 -\kappa_2{\cal H}u^{n+1-k}_i)
 +\dsum_{k=0}^{n}\lambda_k(\kappa_1\delta^2_xu^{n-k}_i
 -\kappa_2{\cal H}u^{n-k}_i)\bigg]\\
 \quad+\tau{\cal H}\phi_i+\dfrac\tau2{\cal H}(f^n_i+f^{n+1}_i),\quad
 1\leq i\leq M-1,\quad 0\leq n\leq N-1.
 \end{array}
 \end{equation}
 The approximate solution is solved with
 \begin{equation}\label{compact-scheme4}
 u_i^0=0,\quad 0\leq i\leq M.
 \end{equation}

 It is easy to see that at each time level, the difference scheme which consists of (\ref{compact-scheme1})--(\ref{compact-scheme4})
 is a linear tridiagonal system with strictly diagonal dominant coefficient matrix.
 Thus the difference scheme has a unique solution.

 \section{Stability and convergence analysis of the compact scheme}
 We give the convergence of the proposed scheme in this section. The main result can be established by the following lemmas:
 \begin{lemma}\label{sequence}{\rm(\cite{Wang})}
 Let $\{\lambda_n\}_{n=0}^\infty$ be defined as {\rm(\ref{sequence-lambda})}, then for any
 positive integer $k$ and real vector
 $(v_1,v_2,\ldots,v_k)^T\in\mathbb{R}^k$,
 it holds that
 $$\dsum_{n=0}^{k-1}\Big(\dsum_{p=0}^n\lambda_pv_{n+1-p}\Big)v_{n+1}\geq0.$$ \end{lemma}
 \begin{lemma}\label{gronwall}{\rm(\cite{Quarteroni})}
 Assume that $\{k_n\}$ and $\{p_n\}$ are nonnegative sequences, and the sequence $\{\phi_n\}$ satisfies
 $$\phi_0\leq g_0,~~~\phi_n\leq g_0+\dsum_{l=0}^{n-1}p_l
 +\dsum_{l=0}^{n-1}k_l\phi_l,~~~n\geq1,$$
 where $g_0\geq0$. Then the sequence $\{\phi_n\}$ satisfies
 $$\phi_n\leq\Big(g_0+\dsum_{l=0}^{n-1}p_l\Big)
 \exp{\Big(\dsum_{l=0}^{n-1}k_l\Big)},~~~n\geq1.$$
 \end{lemma}

 Our compact difference scheme consisting of (\ref{compact-scheme1})--(\ref{compact-scheme4}) has high order convergence.
 To be more precise, we have
 \begin{theorem}\label{con}
 Assume that $u(x,t)\in {\cal C}_{x,t}^{6,2}([0,L]\times[0,T])$ is the solution of {\rm(\ref{main1})--(\ref{main3})} and $\{u_i^k|0\leq i\leq M,~ 0\leq k\leq N\}$
 is a solution of the finite difference scheme {\rm(\ref{compact-scheme1})--(\ref{compact-scheme4})}, respectively. Denote $$e_i^k=u(x_i,t_k)-u_i^k, \quad 0\leq i\leq M, \quad 0\leq k\leq N.$$
 Then there exists a positive constant $c$ such that
 $$\|e^k\|\leq c(\tau^2+h^4), \quad 0\leq k\leq N.$$
 \end{theorem}
 {\bf Proof.} We can easily get the following error equation:
 \begin{equation}\label{error}
 \begin{array}{l}
 \bar C(e^{k+1}-e^{k})=-\dfrac{\kappa_1\tau^{\alpha+1}}{2h^2}\dsum_{l=0}^{k}\lambda_l\bar Q(e^{k+1-l}+e^{k-l})
 -\dfrac{\kappa_2\tau^{\alpha+1}}2\dsum_{l=0}^{k}\lambda_l\bar C(e^{k+1-l}+e^{k-l})
 +\tau \bar R^{k+1},\\
 e_i^0=0, \quad 0\leq i\leq M,
 \end{array}
 \end{equation}
 where $\|\bar R^{k+1}\|\leq c_1(\tau^{2}+h^4)$,
 \begin{equation}\label{matrix1}
 \bar C=\frac1{12}\left(\begin{array}{ccccc}
  10&2&&&\\
  1&10&1&&\\
  &\ddots&\ddots&\ddots&\\
  &&1&10&1\\
  &&&2&10
  \end{array}\right),\quad
  \bar Q=\left(\begin{array}{ccccc}
  2&-2&&&\\
  -1&2&-1&&\\
  &\ddots&\ddots&\ddots&\\
  &&-1&2&-1\\
  &&&-2&2
  \end{array}\right).
  \end{equation}
 Multiplying the equation (\ref{error}) with $\frac12\oplus I\oplus \frac12$, where $I$ is the identity matrix, we get
 \begin{equation}\label{sym-error}
 C(e^{k+1}-e^{k})=-\dfrac{\kappa_1\tau^{\alpha+1}}{2h^2}\dsum_{l=0}^{k}\lambda_l Q(e^{k+1-l}+e^{k-l})
 -\dfrac{\kappa_2\tau^{\alpha+1}}2\dsum_{l=0}^{k}\lambda_lC(e^{k+1-l}+e^{k-l})
 +\tau R^{k+1},
 \end{equation}
 where $\|R^{k+1}\|\leq c_2(\tau^{2}+h^4)$,
 \begin{equation}\label{matrix2}
 C=\dfrac1{12}\left(\begin{array}{ccccc}
  5&1&&&\\
  1&10&1&&\\
   &\ddots&\ddots&\ddots&\\
  &&1&10&1\\
  &&&1&5
 \end{array}\right)=E^2,
 \mbox{ with } E \mbox{ being the square root of } C,
 \end{equation}
 \begin{equation}\label{matrix3}
  Q=\left(\begin{array}{ccccc}
  1&-1&&&\\
  -1&2&-1&&\\
   &\ddots&\ddots&\ddots&\\
  &&-1&2&-1\\
  &&&-1&1
  \end{array}\right)=S^TS,
  \mbox{ with }
  S=\left(\begin{array}{cccc}
  -1&1&&\\
   &\ddots&\ddots\\
  &&-1&1
  \end{array}\right)\in \IR^{M\times{M+1}}.
  \end{equation}
 Here we have used the fact that $C=\frac1{12}tri[1,5,1]+\frac5{12}(0\oplus I\oplus0)$ is positive definite.

 Multiplying (\ref{sym-error}) by $h(e^{k+1}+e^{k})^T$, we obtain
 $$\begin{array}{ll}
 \quad h(e^{k+1}+e^{k})^TC(e^{k+1}-e^{k})\\
 =-\dfrac{\kappa_1\tau^{\alpha+1}}{2h}\dsum_{l=0}^{k}\lambda_l(e^{k+1}+e^{k})^TS^TS(e^{k+1-l}+e^{k-l})
 -\dfrac{\kappa_2h\tau^{\alpha+1}}{2}\dsum_{l=0}^{k}\lambda_l(e^{k+1}+e^{k})^TE^2(e^{k+1-l}+e^{k-l})\\
 \quad+\tau h(e^{k+1}+e^{k})^TR^{k+1}.
 \end{array}$$
 Summing up for $0\leq k\leq n-1$ and noting that $$h(e^{k+1}+e^{k})^TC(e^{k+1}-e^{k})=h\big[({e^{k+1}})^TCe^{k+1}-({e^k})^TCe^{k}\big],
 ~~h({e^n})^TCe^n\geq\dfrac14\|e^n\|^2,$$
 we have, by Lemma \ref{sequence},
 $$\begin{array}{rl}
 \dfrac14\|e^n\|^2&\leq\tau h(e^{n}+e^{n-1})^TR^{n}+\tau h\dsum_{k=0}^{n-2}(e^{k+1}+e^{k})^TR^{k+1}\\
 &\leq\dfrac15\|e^n\|^2+\dfrac{5\tau^2}4\|R^{n}\|^2+\dfrac\tau2\|e^{n-1}\|^2
 +\dfrac{\tau}2\|R^{n}\|^2
 +\dfrac\tau2\dsum_{k=1}^{n-1}\|e^{k}\|^2+\dfrac\tau2\dsum_{k=1}^{n-2}\|e^{k}\|^2
 +\tau\dsum_{k=1}^{n-1}\|R^{k}\|^2\\
 &\leq\dfrac15\|e^n\|^2+\dfrac{5\tau^2}4\|R^{n}\|^2
 +\tau\dsum_{k=1}^{n-1}\|e^{k}\|^2
 +\tau\dsum_{k=1}^{n}\|R^{k}\|^2,
 \end{array}$$
 which gives
 $$\begin{array}{ll}
 \|e^n\|^2&\leq25\tau^2\|R^{n}\|^2
 +20\tau\dsum_{k=1}^{n-1}\|e^{k}\|^2
 +20\tau\dsum_{k=1}^{n}\|R^{k}\|^2\\
 &\leq20\tau\dsum_{k=1}^{n-1}\|e^{k}\|^2+c_3(\tau^2+h^4)^2,
 \end{array}$$
 then the desired result follows by Lemma \ref{gronwall}.$\qquad\Box$
 \begin{remark}\label{re-stability}
 One can adopt the idea of the proof for Theorem {\rm\ref{con}} to
 show that the proposed compact scheme
{\rm(\ref{compact-scheme1})--(\ref{compact-scheme4})} is
unconditionally stable.
 In fact, consider the solution $\{v_i^k\}$
 of
 \begin{equation}\label{sta1}
 \begin{array}{ll}
 \quad{\cal H}(v^{k+1}_0-v^k_0)\\
 =\dfrac{\tau^{\alpha+1}}2
 \bigg[\dsum_{l=0}^{k+1}\lambda_l(\kappa_1\delta^2_xv^{k+1-l}_0
 -\kappa_2{\cal H}v^{k+1-l}_0)
 +\dsum_{l=0}^{k}\lambda_l(\kappa_1\delta^2_xv^{k-l}_0
 -\kappa_2{\cal H}v^{k-l}_0)\bigg]\\
 \quad+\dfrac{\tau^{\alpha+1}}2\dsum_{l=0}^{k+1}
 \lambda_l\bigg(\frac h6(g_x)_0^{k+1-l}-\frac{h^3}{90}
 \Big[(g_{xxx})_0^{k+1-l}+\frac{\kappa_2}{\kappa_1}(g_x)_0^{k+1-l}
 +\frac1{\kappa_1}({_{0}^CD_t^{\alpha+1}g_x})_0^{k+1-l}\Big]\bigg)\\
 \quad+\dfrac{\tau^{\alpha+1}}2\dsum_{l=0}^{k}
 \lambda_l\bigg(\frac h6(g_x)_0^{k-l}-\frac{h^3}{90}
 \Big[(g_{xxx})_0^{k-l}+\frac{\kappa_2}{\kappa_1}(g_x)_0^{k-l}
 +\frac1{\kappa_1}({_{0}^CD_t^{\alpha+1}g_x})_0^{k-l}\Big]\bigg)\\
 \quad+\tau{\cal H}(\phi_0+\tilde\rho_0)+\dfrac\tau2{\cal H}(f^k_0+f^{k+1}_0),\quad 0\leq k\leq N-1,
 \end{array}
 \end{equation}
 \begin{equation}\label{sta3}
 \begin{array}{ll}
 \quad{\cal H}(v^{k+1}_M-v^k_M)\\
 =\dfrac{\tau^{\alpha+1}}2
 \bigg[\dsum_{l=0}^{k+1}\lambda_l(\kappa_1\delta^2_xv^{k+1-l}_M
 -\kappa_2{\cal H}v^{k+1-l}_M)
 +\dsum_{l=0}^{k}\lambda_l(\kappa_1\delta^2_xv^{k-l}_M
 -\kappa_2{\cal H}v^{k-l}_M)\bigg]\\
 \quad-\dfrac{\tau^{\alpha+1}}2\dsum_{l=0}^{k+1}
 \lambda_l\bigg(\frac h6(g_x)_M^{k+1-l}-\frac{h^3}{90}
 \Big[(g_{xxx})_M^{k+1-l}+\frac{\kappa_2}{\kappa_1}(g_x)_M^{k+1-l}
 +\frac1{\kappa_1}({_{0}^CD_t^{\alpha+1}g_x})_M^{k+1-l}\Big]\bigg)\\
 \quad-\dfrac{\tau^{\alpha+1}}2\dsum_{l=0}^{k}
 \lambda_l\bigg(\frac h6(g_x)_M^{k-l}-\frac{h^3}{90}
 \Big[(g_{xxx})_M^{k-l}+\frac{\kappa_2}{\kappa_1}(g_x)_M^{k-l}
 +\frac1{\kappa_1}({_{0}^CD_t^{\alpha+1}g_x})_M^{k-l}\Big]\bigg)\\
 \quad+\tau{\cal H}(\phi_M+\tilde\rho_M)+\dfrac\tau2{\cal H}(f^k_M+f^{k+1}_M),\quad 0\leq k\leq N-1,
 \end{array}
 \end{equation}
 \begin{equation}\label{sta2}
 \begin{array}{ll}
 \quad{\cal H}(v^{k+1}_i-v^k_i)\\
 =\dfrac{\tau^{\alpha+1}}2
 \bigg[\dsum_{l=0}^{k+1}\lambda_l(\kappa_1\delta^2_xv^{k+1-l}_i
 -\kappa_2{\cal H}v^{k+1-l}_i)
 +\dsum_{l=0}^{k}\lambda_l(\kappa_1\delta^2_xv^{k-l}_i
 -\kappa_2{\cal H}v^{k-l}_i)\bigg]\qquad\qquad\quad\\
 \quad+\tau{\cal H}(\phi_i+\tilde\rho_i)+\dfrac\tau2{\cal H}(f^k_i+f^{k+1}_i),\quad
 1\leq i\leq M-1,\quad 0\leq k\leq N-1,
 \end{array}
 \end{equation}
 with $v_i^0=\rho_i,~0\leq i\leq M$.
 Then, by  {\rm(\ref{compact-scheme1})--(\ref{compact-scheme2})}
 and {\rm(\ref{sta1})--(\ref{sta2})},
 one can check that $\varepsilon_i^l=v_i^l-u_i^l-\rho_i$ satisfy
 \begin{equation}\label{sta5}
 \begin{array}{l}
 {\cal H}(\varepsilon^{k+1}_i-\varepsilon^k_i)=\dfrac{\tau^{\alpha+1}}2
 \bigg[\dsum_{l=0}^{k+1}\lambda_l(\kappa_1\delta^2_x\varepsilon^{k+1-l}_i
 -\kappa_2{\cal H}\varepsilon^{k+1-l}_i)
 +\dsum_{l=0}^{k}\lambda_l(\kappa_1\delta^2_x\varepsilon^{k-l}_i
 -\kappa_2{\cal H}\varepsilon^{k-l}_i)\bigg]+\tau{\cal H}\tilde\rho_i\\
 ~+\dfrac{\tau^{\alpha+1}}2\bigg[\dsum_{l=0}^{k+1}\lambda_l
 (\kappa_1\delta^2_x\rho_i-\kappa_2{\cal H}\rho_i)
 +\dsum_{l=0}^{k}\lambda_l(\kappa_1\delta^2_x\rho_i-\kappa_2{\cal H}\rho_i)\bigg],
 \quad 0\leq i\leq M,~~0\leq k\leq N-1,\\
 \varepsilon_i^0=0,\quad 0\leq i\leq M.
 \end{array}
 \end{equation}
 By following the proof for Theorem {\rm\ref{con}}
 and noting $\tau^\alpha\sum\limits_{l=0}^{k+1}\lambda_l=\frac1{\Gamma(\alpha+1)}+O(\tau)$,
 we then have the estimate
 \begin{equation}\label{stability}
 \begin{array}{ll}
 \|\varepsilon^k\|^2&\leq20\tau\dsum_{l=0}^{k-1}\|\varepsilon^{l}\|^2
 +\Big[\frac5{\Gamma(\alpha+1)}+1\Big]^2\Big[\|\kappa_1\delta_x^2\rho\|^2
 +\|\kappa_2\rho\|^2+\|\tilde\rho\|^2\Big]\\
 &\leq e^{20T}\Big[\dfrac5{\Gamma(\alpha+1)}+1\Big]^2
 \Big[\|\kappa_1\delta_x^2\rho\|^2
 +\|\kappa_2\rho\|^2+\|\tilde\rho\|^2\Big].
 \end{array}
 \end{equation}
 This implies
 \begin{align*}
 \|v^k-u^k\|&\leq\|v^k-u^k-\rho\|+\|\rho\|\\
 &\leq e^{10T}\Big[\dfrac5{\Gamma(\alpha+1)}+1\Big]
 \sqrt{\|\kappa_1\delta_x^2\rho\|^2+\|\kappa_2\rho\|^2+\|\tilde\rho\|^2}+\|\rho\|,
 \end{align*}
 concluding the stability of the scheme.
 \end{remark}

 \section{The compact ADI scheme for the two-dimensional problem}
 In this section, we turn to study the two-dimensional problem:
 \begin{align}\label{2Dmain1}
 &_{0}^CD_t^{\gamma}u
 =\Delta u-u+g(x,y,t), \quad (x,y)\in \Omega, \quad 0<t\leq T, \quad 1<\gamma<2,\\\label{2Dmain2}
 &u(x,y,0)=0, \quad \frac{\partial u(x,y,0)}{\partial t}=\phi(x,y), \quad (x,y)\in\bar{\Omega}=\Omega\cup\partial\Omega,\\\label{2Dmain3}
 &\dfrac{\partial u(x,y,t)}{\partial n}\Big|_{\partial\Omega}=0, \quad (x,y)\in \partial\Omega, \quad 0<t\leq T,
 \end{align}
 where $\Delta$ is the
 two-dimensional Laplacian, $n$ is the unit outward normal vector of the domain $\Omega=(0,L_1)\times(0,L_2)$
 with boundary $\partial\Omega$.

 An equivalent form of (\ref{2Dmain1}) read as:
 \begin{equation}\label{2Dchanged-main}
 \dfrac{\partial u(x,y,t)}{\partial t}=\phi(x,y)+
 \dfrac1{\Gamma(\alpha)}\dint_0^t(t-s)^{\alpha-1}
 [\Delta u(x,y,s)-u(x,y,s)]ds+f(x,y,t),
 \end{equation}
 where $(x,y)\in\Omega,~0<t\leq T,~0<\alpha=\gamma-1<1,~f(x,y,t)={_{0}}I^\alpha_tg(x,y,t)$.

 Discretization of (\ref{2Dchanged-main}) are carried out with steps similar to that of the one-dimensional problem.
 To this end, we let $h_1=\frac{L_1}{M_1},~h_2=\frac{L_2}{M_2}$ and $\tau=\frac TN$ be the spatial and temporal step sizes respectively, where  $M_1,~M_2$ and $N$ are some given integers.
 For $i=0,1,\ldots,M_1,~j=0,1,\ldots,M_2$ and $k=0,1,\ldots,N$, denote $x_i=ih_1,~y_j=jh_2,~t_k=k\tau$.
 We introduce the following notations
 on a grid function $u=\{u_{ij}^k|0\leq i\leq M_1,~0\leq j\leq M_2,~0\leq k\leq N\}$:
 $$\begin{array}{c}
 \delta_xu_{i-\frac12,j}=\frac1{h_1}(u_{ij}-u_{i-1,j}),
 \end{array}$$
 $$\delta_x^2u_{ij}=\left\{
 \begin{array}{ll}
 \frac2{h_1}\delta_xu_{\frac12,j},&i=0,~~0\leq j\leq M_2,\\
 \frac1{h_1}(\delta_xu_{i+\frac12,j}-\delta_xu_{i-\frac12,j}),&1\leq i\leq M_1-1,~~0\leq j\leq M_2,\\
 -\frac2{h_1}\delta_xu_{M_1-\frac12,j},&i=M_1,~~0\leq j\leq M_2,
 \end{array}\right.$$
 $${\cal H}_xu_{ij}=\left\{
 \begin{array}{ll}
 \frac16(5u_{0,j}+u_{1,j}), &i=0,~~0\leq j\leq M_2,\\
 \frac1{12}(u_{i-1,j}+10u_{i,j}+u_{i+1,j}), &1\leq i\leq M_1-1,~~0\leq j\leq M_2,\\
 \frac16(u_{M_1-1,j}+5u_{M_1,j}), &i=M_1,~~0\leq j\leq M_2.
 \end{array}\right.$$
 One can defined similar notations in the $y$ direction. We further denote:
 $$\begin{array}{c}
 {\cal H}u_{ij}={\cal H}_x{\cal H}_yu_{ij},~~\Lambda u_{ij}=({\cal H}_y\delta^2_x+{\cal H}_x\delta_y^2)u_{ij},\\
 \langle u,v\rangle=h_1h_2\sum\limits_{i=0}^{M_1}\sum\limits_{j=0}^{M_2}u_{ij}v_{ij},~~\|u\|^2=\langle u,u\rangle,~~\|u\|_\infty=\max\limits_{0\leq i\leq M_1,~0\leq j\leq M_2}|u_{ij}|.
 \end{array}$$
 With all the preparation, we now  give the compact ADI scheme.

 We first denote $\mu=\frac{\tau^{\alpha+1}}2$ and
 $G_{ij}^n=(G_{ij}^1)^n+(G_{ij}^2)^n,$ where
 $$(G_{ij}^1)^n=\left\{
 \begin{array}{ll}
 \mu\dsum_{k=0}^{n}
 \lambda_k\bigg(\frac{h_1}6{\cal H}_y(g_x)_{0,j}^{n-k}
 -\frac{h_1^3}{90}{\cal H}_y\Big[(g_{xxx})_{0,j}^{n-k}-(g_{xyy})_{0,j}^{n-k}
 +(g_{x})_{0,j}^{n-k}+({_{0}^CD_t^{\alpha+1}g_x})_{0,j}^{n-k}\Big]\bigg),\\
 ~i=0,\quad 0\leq j\leq M_2,\\[5pt]
 0,\quad 1\leq i\leq M_1-1,\quad 0\leq j\leq M_2,\\[5pt]
 -\mu\dsum_{k=0}^{n}
 \lambda_k\bigg(\frac{h_1}6{\cal H}_y(g_x)_{M_1,j}^{n-k}
 -\frac{h_1^3}{90}{\cal H}_y\Big[(g_{xxx})_{M_1,j}^{n-k}-(g_{xyy})_{M_1,j}^{n-k}
 +(g_{x})_{M_1,j}^{n-k}+({_{0}^CD_t^{\alpha+1}g_x})_{M_1,j}^{n-k}\Big]\bigg),\\
 ~i=M_1,\quad 0\leq j\leq M_2,
 \end{array}\right.$$ 
 $$(G_{ij}^2)^n=\left\{
 \begin{array}{ll}
 \mu\dsum_{k=0}^{n}
 \lambda_k\bigg(\frac{h_2}6{\cal H}_x(g_y)_{i,0}^{n-k}
 -\frac{h_2^3}{90}{\cal H}_x\Big[(g_{yyy})_{i,0}^{n-k}-(g_{xxy})_{i,0}^{n-k}
 +(g_{y})_{i,0}^{n-k}+({_{0}^CD_t^{\alpha+1}g_y})_{i,0}^{n-k}\Big]\bigg),\\
 ~j=0,\quad 0\leq i\leq M_1,\\[5pt]
 0,\quad 1\leq j\leq M_2-1,\quad 0\leq i\leq M_1,\\[5pt]
 -\mu\dsum_{k=0}^{n}
 \lambda_k\bigg(\frac{h_2}6{\cal H}_x(g_y)_{i,M_2}^{n-k}
 -\frac{h_2^3}{90}{\cal H}_x\Big[(g_{yyy})_{i,M_2}^{n-k}-(g_{xxy})_{i,M_2}^{n-k}
 +(g_{y})_{i,M_2}^{n-k}+({_{0}^CD_t^{\alpha+1}g_y})_{i,M_2}^{n-k}\Big]\bigg),\\
 ~j=M_2,\quad 0\leq i\leq M_1.
 \end{array}\right.$$
 Following the steps in  the one dimensional case, one can deduce the following:
 \begin{equation}\label{2Dcompact1}
 \begin{array}{ll}
 {\cal H}(u^{n+1}_{ij}-u^n_{ij})=\tau{\cal H}\phi_{ij}
 +\mu\bigg[\dsum_{k=0}^{n+1}\lambda_k(\Lambda-{\cal H})u^{n+1-k}_{ij}
 +\dsum_{k=0}^{n}\lambda_k(\Lambda-{\cal H})u^{n-k}_{ij}\bigg]\\
 \qquad\qquad\qquad\quad~+\dfrac\tau2{\cal H}(f^n_{ij}+f^{n+1}_{ij})
 +\dfrac12(G_{ij}^n+G_{ij}^{n+1})+\tau(R_1)_{ij}^{n+1},\\
 u^0_{ij}=0,\quad (x_i,y_j)\in\bar{\Omega},
 \end{array}
 \end{equation}
 where $(R_1)_{ij}^{n+1}=O(\tau^2+h_1^4+h_2^4)$.

 \medskip
 Denoting $F_{ij}^n=\frac12(\tau{\cal H}f_{ij}^n+G_{ij}^n)$, and adding a small term $\frac{\mu^2\lambda_0^2}{1+\mu\lambda_0}\delta_x^2\delta_y^2(u_{ij}^{n+1}-u_{ij}^{n})=O(\tau^{3+2\alpha})$
 on both sides of (\ref{2Dcompact1}), we have
 \begin{equation}\label{2Dcompact-scheme}
 \begin{array}{ll}
 \quad{\cal H}(u^{n+1}_{ij}-u^n_{ij})+
 \frac{\mu^2\lambda_0^2}{1+\mu\lambda_0}\delta_x^2\delta_y^2(u_{ij}^{n+1}-u_{ij}^{n})\\
 =\tau{\cal H}\phi_{ij}
 +\mu\bigg[\dsum_{k=0}^{n+1}\lambda_k(\Lambda-{\cal H})u^{n+1-k}_{ij}
 +\dsum_{k=0}^{n}\lambda_k(\Lambda-{\cal H})u^{n-k}_{ij}\bigg]
 +F^n_{ij}+F^{n+1}_{ij}+\tau R_{ij}^{n+1}\\
 u^0_{ij}=0,\quad (x_i,y_j)\in\bar{\Omega},
 \end{array}
 \end{equation}
 with $R_{ij}^{n+1}=O(\tau^{2}+h_1^4+h_2^4)$. Omitting the truncation error in (\ref{2Dcompact-scheme}), we reach the following
 scheme in the ADI setting:
 $$\begin{array}{ll}
 \quad\big(\sqrt{1+\mu\lambda_0}{\cal H}_x-\frac{\mu\lambda_0}{\sqrt{1+\mu\lambda_0}}\delta^2_x\big)
 \big(\sqrt{1+\mu\lambda_0}{\cal H}_y-\frac{\mu\lambda_0}{\sqrt{1+\mu\lambda_0}}\delta_y^2\big)u_{ij}^{n+1}\\
 ={\cal H}u_{ij}^{n}+\frac{\mu^2\lambda_0^2}{1+\mu\lambda_0}\delta_x^2\delta_y^2u_{ij}^{n}
 +\mu\bigg[\dsum_{k=1}^{n+1}\lambda_k(\Lambda-{\cal H})u^{n+1-k}_{ij}
 +\dsum_{k=0}^{n}\lambda_k(\Lambda-{\cal H})u^{n-k}_{ij}\bigg]
 +\tau{\cal H}\phi_{ij}+F^n_{ij}+F^{n+1}_{ij},\\
 \quad(x_i,y_j)\in\Omega,~~0\leq n\leq N-1,\\
 u_{ij}^0=0, \quad (x_i,y_j)\in\bar{\Omega}.
 \end{array}$$
 For ADI methods (see \cite{Zhang2} for example),
 the solution $\{u_{ij}^{n+1}\}$ is determined by solving two independent one-dimensional problems.
 Specifically,  the intermediate variables
 $$\begin{array}{rl}
 u_{ij}^*=\big(\sqrt{1+\mu\lambda_0}{\cal H}_y-\frac{\mu\lambda_0}{\sqrt{1+\mu\lambda_0}}\delta_y^2\big)u_{ij}^{n+1},~~0\leq i\leq M_1,~0\leq j\leq M_2,
 \end{array}$$
 are first solved from the following system with fixed $j\in\{0,1,\ldots,M_2\}$:
 $$\begin{array}{ll}
 \quad\big(\sqrt{1+\mu\lambda_0}{\cal H}_x-\frac{\mu\lambda_0}{\sqrt{1+\mu\lambda_0}}\delta^2_x\big)u_{ij}^*\\
 ={\cal H}u_{ij}^{n}+\frac{\mu^2\lambda_0^2}{1+\mu\lambda_0}\delta_x^2\delta_y^2u_{ij}^{n}
 +\mu\bigg[\dsum_{k=1}^{n+1}\lambda_k(\Lambda-{\cal H})u^{n+1-k}_{ij}
 +\dsum_{k=0}^{n}\lambda_k(\Lambda-{\cal H})u^{n-k}_{ij}\bigg]+\tau{\cal H}\phi_{ij}+F^n_{ij}+F^{n+1}_{ij},\\
 \quad0\leq i\leq M_1,
 \end{array}$$

 When $\{u_{ij}^*\}$ is ready, the approximate solution $\{u_{ij}^{n+1}\}$ is solved from the following system for fixed $i\in\{0,1,\ldots,M_1\}$:
 $$\Big(\sqrt{1+\mu\lambda_0}{\cal H}_y-\frac{\mu\lambda_0}{\sqrt{1+\mu\lambda_0}}\delta_y^2\Big)u_{ij}^{n+1}=u_{ij}^*,
 ~~0\leq j\leq M_2.$$
 By implementing the ADI method, the computational cost for solving a two-dimensional problem can be greatly reduced.

 \medskip
 We now proceed to give the convergence result of our compact ADI scheme (\ref{2Dcompact-scheme}).
 We remark that, with the convergence of the scheme,
 one can show that the scheme is stable in the same sense as that given in Remark \ref{re-stability}.
 \begin{theorem}
 Assume that $u(x,y,t)\in {\cal C}_{x,y,t}^{6,6,2}(\Omega\times[0,T])$ is the solution of {\rm(\ref{2Dmain1})--(\ref{2Dmain3})} and $\{u_{ij}^k|0\leq i\leq M_1,~ 0\leq j\leq M_2,~ 0\leq k\leq N\}$
 is a solution of the finite difference scheme {\rm(\ref{2Dcompact-scheme})}, respectively. Denote $$e_{ij}^k=u(x_i,y_j,t_k)-u_{ij}^k, \quad 0\leq i\leq M_1,~0\leq j\leq M_2,~1\leq k\leq N.$$
 Then there exists a positive constant $\tilde c$ such that
 $$\|e^k\|\leq\tilde c(\tau^2+h_1^4+h_2^4), \quad 0\leq k\leq N,$$
 where $e^k=[e^k_{0,0},e^k_{1,0},\cdots,e^k_{M_1,0},e^k_{0,1},e^k_{1,1},\cdots,e^k_{M_1,1},\cdots,e^k_{0,M_2},e^k_{1,M_2},\cdots,e^k_{M_1,M_2}]^T$.
 \end{theorem}
 {\bf Proof.} One can easily check that the following error equation holds:
 \begin{equation}\label{2Derror}
 \begin{array}{rl}
 &\quad(\bar C_{M_2+1}\otimes\bar C_{M_1+1})(e^{k+1}-e^{k})+
 \frac{\mu^2\lambda_0^2}{(1+\mu\lambda_0)h_1^2h_2^2}(\bar Q_{M_2+1}\otimes\bar Q_{M_1+1})(e^{k+1}-e^{k})\\
 &=-\dfrac{\tau^{\alpha+1}}{2h_1^2}\dsum_{l=0}^{k}\lambda_l(\bar C_{M_2+1}\otimes\bar Q_{M_1+1})(e^{k+1-l}+e^{k-l})
 -\dfrac{\tau^{\alpha+1}}{2h_2^2}\dsum_{l=0}^{k}\lambda_l(\bar Q_{M_2+1}\otimes\bar C_{M_1+1})(e^{k+1-l}+e^{k-l})\\
 &\quad-\dfrac{\tau^{\alpha+1}}2\dsum_{l=0}^{k}\lambda_l(\bar C_{M_2+1}\otimes\bar C_{M_1+1})(e^{k+1-l}+e^{k-l})
 +\tau \bar R^{k+1},\\
 &e_{ij}^0=0, \quad 0\leq i\leq M_1,\quad 0\leq j\leq M_2,
 \end{array}
 \end{equation}
 where $\|\bar R^{k+1}\|\leq\tilde c_1(\tau^{2}+h_1^4+h_2^4)$,
 and the matrices $\bar C_{M_1+1},~\bar C_{M_2+1},~\bar Q_{M_1+1},~\bar Q_{M_2+1}$ are given in (\ref{matrix1})
 with the corresponding sizes given by the subscripts.

 \medskip
 Multiplying the equation (\ref{2Derror}) with $(\frac12\oplus I_{M_2-1}\oplus \frac12)\otimes(\frac12\oplus I_{M_1-1}\oplus \frac12)$, we get
 \begin{equation}\label{2Dsym-error}
 \begin{array}{rl}
 &\quad(C_{M_2+1}\otimes C_{M_1+1})(e^{k+1}-e^{k})+
 \frac{\mu^2\lambda_0^2}{(1+\mu\lambda_0)h_1^2h_2^2}(Q_{M_2+1}\otimes Q_{M_1+1})(e^{k+1}-e^{k})\\
 &=-\dfrac{\tau^{\alpha+1}}{2h_1^2}\dsum_{l=0}^{k}\lambda_l(C_{M_2+1}\otimes Q_{M_1+1})(e^{k+1-l}+e^{k-l})-\dfrac{\tau^{\alpha+1}}{2h_2^2}\dsum_{l=0}^{k}\lambda_l(Q_{M_2+1}\otimes C_{M_1+1})(e^{k+1-l}+e^{k-l})\\
 &\quad-\dfrac{\tau^{\alpha+1}}2\dsum_{l=0}^{k}\lambda_l(C_{M_2+1}\otimes C_{M_1+1})(e^{k+1-l}+e^{k-l})
 +\tau R^{k+1},\\
 &e_{ij}^0=0, \quad 0\leq i\leq M_1, \quad 0\leq j\leq M_2,
 \end{array}
 \end{equation}
 where $\|R^{k+1}\|\leq\tilde c_2(\tau^{2}+h_1^4+h_2^4)$
 and $$C_{M_1+1}=E_{M_1+1}^2,~C_{M_2+1}=E_{M_2+1}^2,~Q_{M_1+1}=S_{M_1+1}^TS_{M_1+1},
 ~Q_{M_2+1}=S_{M_2+1}^TS_{M_2+1}$$
 are as in (\ref{matrix2}), (\ref{matrix3}) respectively.
 Once again, we have used the subscripts to indicate the sizes of the matrices.

 We can now multiply (\ref{2Dsym-error}) by $h_1h_2(e^{k+1}+e^{k})^T$ and
 add up the equations for $0\leq k\leq n-1$.
 Noting that $$(e^{k+1}+e^{k})^T(C_{M_2+1}\otimes C_{M_1+1})(e^{k+1}-e^{k})
 =({e^{k+1}})^T(C_{M_2+1}\otimes C_{M_1+1})e^{k+1}-({e^k})^T(C_{M_2+1}\otimes C_{M_1+1})e^{k},$$
 $$(e^{k+1}+e^{k})^T(Q_{M_2+1}\otimes Q_{M_1+1})(e^{k+1}-e^{k})
 =({e^{k+1}})^T(Q_{M_2+1}\otimes Q_{M_1+1})e^{k+1}-({e^k})^T(Q_{M_2+1}\otimes Q_{M_1+1})e^{k},$$
 $$h_1h_2({e^n})^T(C_{M_2+1}\otimes C_{M_1+1})e^n\geq\dfrac1{16}\|e^n\|^2,
 \quad ({e^n})^T(Q_{M_2+1}\otimes Q_{M_1+1})e^n\geq0,$$
 $$(C_{M_2+1}\otimes Q_{M_1+1})=(E_{M_2+1}\otimes S_{M_1+1}^T)(E_{M_2+1}\otimes S_{M_1+1}),$$
 $$(Q_{M_2+1}\otimes C_{M_1+1})=(S_{M_2+1}^T\otimes E_{M_1+1})(S_{M_2+1}\otimes E_{M_1+1}),$$
 we have, by Lemma \ref{sequence}, that
 $$\begin{array}{rl}
 \dfrac1{16}\|e^n\|^2&\leq\tau h_1h_2(e^{n}+e^{n-1})^TR^{n}+\tau h_1h_2\dsum_{k=0}^{n-2}(e^{k+1}+e^{k})^TR^{k+1}\\
 &\leq\dfrac1{18}\|e^n\|^2+\dfrac{9\tau^2}2\|R^{n}\|^2+\dfrac\tau2\|e^{n-1}\|^2
 +\dfrac{\tau}2\|R^{n}\|^2
 +\dfrac\tau2\dsum_{k=1}^{n-1}\|e^{k}\|^2+\dfrac\tau2\dsum_{k=1}^{n-2}\|e^{k}\|^2
 +\tau\dsum_{k=1}^{n-1}\|R^{k}\|^2\\
 &\leq\dfrac1{18}\|e^n\|^2+\dfrac{9\tau^2}2\|R^{n}\|^2
 +\tau\dsum_{k=1}^{n-1}\|e^{k}\|^2
 +\tau\dsum_{k=1}^{n}\|R^{k}\|^2,
 \end{array}$$
 from which we can conclude the theorem just as in the one dimensional case. $\qquad\Box$

 \section{Numerical experiments}
 In this section, we carry out numerical experiments for the finite difference scheme to illustrate our theoretical statements. All our tests were done in MATLAB.
 Although our theoretical results are given by the discrete $L^2$ norm, we find that the maximum norm errors
 $$E_\infty(h,\tau)=\max\limits_{0\leq k\leq N}\|U^k-u^k\|_\infty$$
 between the exact and the numerical solutions also match the proposed order for the examples we have tested
 (we remark here that we have similar observations in \cite{Wang}).
 Therefore, in the numerical examples given below, the maximum norm errors are reported.

 \medskip
 We first consider the following one-dimensional problem:
 \begin{example}\label{ex1}
 $$\begin{array}{rl}
 &_{0}^CD_t^{\gamma}u
 =\frac{\partial^2u}{\partial x^2}-u+g(x,t),
 \quad 0\leq x\leq 1, \quad 0<t\leq 1, \quad 1<\gamma<2,\\
 &u(x,0)=0, \quad \frac{\partial u(x,0)}{\partial t}=0, \quad 0\leq x\leq 1,\\
 &\frac{\partial u(0,t)}{\partial x}=0, \quad \frac{\partial u(L,t)}{\partial x}=0, \quad 0<t\leq 1,
 \end{array}$$
 where
 $g(x,t)=\frac{\Gamma(\gamma+3)}2t^2e^xx^2(1-x)^2-e^xt^{\gamma+2}(2-8x+8x^3)$.

 Note that the equation can be equivalently written as
 $$\begin{array}{rl}
 \frac{\partial u(x,t)}{\partial t}=
 {_{0}I^\alpha_t[u_{xx}(x,t)-u(x,t)]}+f(x,t), \quad 0\leq x\leq 1, \quad 0<t\leq 1,
 \end{array}$$
 where $\alpha=\gamma-1,~f(x,t)=(\alpha+3)e^xx^2(1-x)^2t^{\alpha+2}
 -\frac{\Gamma(\alpha+4)}{\Gamma(2\alpha+4)}e^x(2-8x+8x^3)t^{2\alpha+3}$.
 The exact solution is $u(x,t)=e^xx^2(1-x)^2t^{\alpha+3}$.
 \end{example}

 Figure 1 plots the curves of the exact solution and numerical solution for the problem at $t=1$ with $\alpha=0.5,~h=\tau=\frac1{100}$.
 The maximum norm errors
 are shown in Table \ref{table1} and Table \ref{table2}.
 Furthermore, the temporal convergence order and spatial convergence order, denoted by
 $$Rate1=\log_2\bigg(\frac{E_\infty(h,2\tau)}{E_\infty(h,\tau)}\bigg)~~\mbox{and}~~ Rate2=\log_2\bigg(\frac{E_\infty(2h,\tau)}{E_\infty(h,\tau)}\bigg),$$
 respectively, are reported.

 \begin{figure}\label{figure1}
 \begin{center}
 \includegraphics[height=7cm]{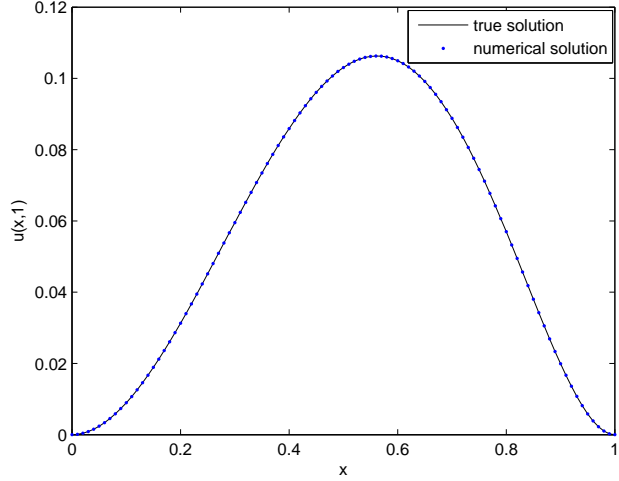}
 \end{center}
 \caption{The exact solution and numerical solution for Example \ref{ex1} at $t=1$, when $\alpha=0.5,~h=\tau=\frac1{100}$.}
 \end{figure}

 \begin{table}[hbt!]
\begin{center}
\caption{Numerical convergence orders in temporal direction with $h=\frac1{50}$ for Example \ref{ex1}.}\label{table1}
\renewcommand{\arraystretch}{0.8}
\def\temptablewidth{0.85\textwidth}
{\rule{\temptablewidth}{0.7pt}}
\begin{tabular*}{\temptablewidth}{@{\extracolsep{\fill}}lllllll}
 $\tau$     &\multicolumn{2}{c}{$\alpha=0.3$}&\multicolumn{2}{c}{$\alpha=0.5$}
&\multicolumn{2}{c}{$\alpha=0.7$}\\
 \cline{2-3}\cline{4-5}\cline{6-7}
&$E_{\infty}(h,\tau)$&Rate1&$E_{\infty}(h,\tau)$
&Rate1&$E_{\infty}(h,\tau)$&Rate1\\\hline
$1/5$   & 1.6417e-3  & $\ast$  & 2.3844e-3  & $\ast$  &3.1904e-3 & $\ast$\\
$1/10$  & 4.1558e-4  & 1.9820  & 6.0481e-4  & 1.9791  &7.9961e-4 & 1.9964\\
$1/20$  & 1.0441e-4  & 1.9929  & 1.5221e-4  & 1.9904  &2.0066e-4 & 1.9945\\
$1/40$  & 2.6115e-5  & 1.9993  & 3.8133e-5  & 1.9970  &5.0254e-5 & 1.9975\\
$1/80$  & 6.4822e-6  & 2.0103  & 9.5080e-6  & 2.0038  &1.2551e-5 & 2.0015
\end{tabular*}
{\rule{\temptablewidth}{0.7pt}}
\end{center}
\end{table}

 \begin{table}[hbt!]
 \begin{center}
 \caption{Numerical convergence orders in spatial direction with $\tau=\frac1{2000}$ when $\alpha=0.5$ for Example \ref{ex1}.}
 \label{table2}
 \renewcommand{\arraystretch}{0.8}
 \def\temptablewidth{0.5\textwidth}
{\rule{\temptablewidth}{0.7pt}}
 \begin{tabular*}{\temptablewidth}{@{\extracolsep{\fill}}lll}
 $h$  & $E_\infty(h,\tau)$& Rate2\\\hline
 1/2  & 2.2688e-2    & $\ast$\\
 1/4  & 1.3235e-3    & 4.0995\\
 1/8  & 8.2429e-5    & 4.0050\\
 1/16 & 5.1310e-6    & 4.0058\\
 1/32 & 3.0926e-7    & 4.0524
 \end{tabular*}
 {\rule{\temptablewidth}{0.7pt}}
 \end{center}
 \end{table}

 Next we turn to consider the stability of the scheme
 by testing (\ref{stability}) numerically.
 We note that the bound in (\ref{stability}) has been magnified to a
 certain extend when it is derived theoretically.
 In our test, we find that the mere quantity $B\dot=\big[\frac5{\Gamma(\alpha+1)}+1\big]
 \sqrt{\|\delta_x^2\rho\|^2+\|\rho\|^2+\|\tilde\rho\|^2}$
 already serves as a good bound for $\|\varepsilon^N\|$.
 We have considered two kinds of perturbation given by the discretization of some functions
 $\rho$, $\tilde\rho$
 and the results are given in Table \ref{table0}.
 \begin{table}[hbt!]
 \begin{center}
 \caption{Stability of the scheme for Example \ref{ex1} when $T=1,~\alpha=0.5$.}
 \label{table0}
 \renewcommand{\arraystretch}{0.99}
 \def\temptablewidth{0.8\textwidth}
 {\rule{\temptablewidth}{0.7pt}}
 \begin{tabular*}{\temptablewidth}{@{\extracolsep{\fill}}lllllllll}
 \multicolumn{4}{c}{$\rho=\tilde\rho=0.1x$} & & \multicolumn{4}{c}{$\rho=\tilde\rho=0.1\sin(x)$}\\\hline
 $M=\frac1h$ & $N=\frac1{\tau}$ & $\|\varepsilon^N\|$ & $B$&
 &$M=\frac1h$ & $N=\frac1{\tau}$ & $\|\varepsilon^N\|$ & $B$\\\hline
 $100$&$500$ &0.3097 &5.4637 & & $100$& $500$ &0.2668 & 6.0461\\
       & $1000$ &0.3102 & 5.4637 && & $1000$ &0.2672 & 6.0461\\
       & $2000$ &0.3105 & 5.4637 && & $2000$ &0.2674 & 6.0461\\
       & $4000$ &0.3107 & 5.4637 && & $4000$ &0.2675 & 6.0461\\\hline
 $200$&$500$ &0.4349 &7.6982 & & $200$& $500$ &0.3747 & 8.5231\\
       & $1000$ &0.4356 & 7.6982 && & $1000$ &0.3752 & 8.5231\\
       & $2000$ &0.4360 & 7.6982 && & $2000$ &0.3755 & 8.5231\\
       & $4000$ &0.4363 & 7.6982 && & $4000$ &0.3757 & 8.5231
 \end{tabular*}
 {\rule{\temptablewidth}{0.7pt}}
 \end{center}
 \end{table}

 \bigskip
 The next example is a two-dimensional problem.
 \begin{example}\label{ex2}
 $$\begin{array}{rl}
 &_{0}^CD_t^{\gamma}u
 =\Delta u-u
 +\cos(x)\cos(y)\big[\frac{\Gamma(\gamma+4)}6t^3+3t^{\gamma+3}\big],
 ~~(x,y)\in\Omega=(0,\pi)\times(0,\pi),~~0<t\leq 1,\\[1pt]
 &u(x,y,0)=0,~~\frac{\partial u(x,y,0)}{\partial t}=0,~~ (x,y)\in\bar{\Omega},\\[1pt]
 &\frac{\partial u(x,y,t)}{\partial n}\big|_{\partial\Omega}=0,~~(x,y)\in \partial\Omega,~~0<t\leq1.
 \end{array}$$
 Note that the equation can be equivalently written as
 $$\begin{array}{rl}
 \frac{\partial u(x,y,t)}{\partial t}=
 {_{0}I^\alpha_t(\Delta u-u)}+\cos(x)\cos(y)\big[(\alpha+4)t^{\alpha+3}
 +\frac{3\Gamma(\alpha+5)}{\Gamma(2\alpha+5)}t^{2\alpha+4}\big],
 \end{array}$$
 where $\alpha=\gamma-1$. The exact solution for this problem is $u(x,t)=\cos(x)\cos(y)t^{\alpha+4}$.
 \end{example}

 We let $h_1=h_2=h$, in this example. Figure 2 shows the exact solution (left) and numerical solution (right) for Example \ref{ex2}, when $\alpha=0.5,~h=\tau=\frac1{50}$. In addition, the maximum norm errors between the exact and the numerical solutions
 $$E_\infty(h,\tau)=\max\limits_{0\leq k\leq N}
 \max\limits_{(x_i,y_j)\in\Omega}|u(x_i,y_j,t_k)-u_{ij}^k|$$
 are shown in Table \ref{table3} and Table \ref{table4}. Meanwhile, the temporal convergence order and spatial convergence order, denoted by
 $$Rate1=\log_2\bigg(\frac{E_\infty(h,2\tau)}{E_\infty(h,\tau)}\bigg)
 ~~\mbox{and}~~ Rate2=\log_2\bigg(\frac{E_\infty(2h,\tau)}{E_\infty(h,\tau)}\bigg),$$
 respectively, are reported. These tables confirm the theoretical analysis.

 \begin{figure}\label{figure2}
 \begin{center}
 \includegraphics[height=6.3cm,width=17.3cm]{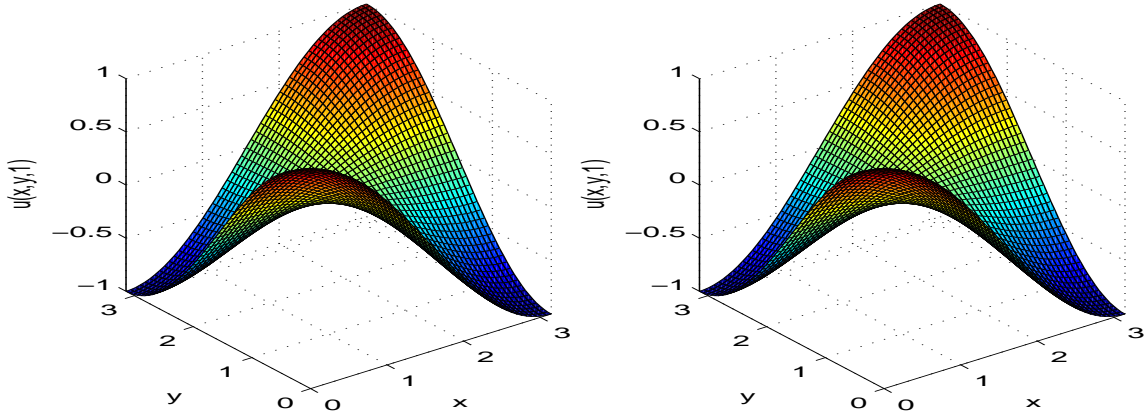}
 \end{center}
 \caption{The exact solution (left) and numerical solution (right) for Example \ref{ex2}, when $\alpha=0.5,~h=\tau=\frac1{50}$.}
 \end{figure}

 \begin{table}[hbt!]
\begin{center}
\caption{Numerical convergence orders in temporal direction with $h=\frac{\pi}{50}$ for Example \ref{ex2}.}\label{table3}
\renewcommand{\arraystretch}{0.8}
\def\temptablewidth{0.85\textwidth}
{\rule{\temptablewidth}{0.7pt}}
\begin{tabular*}{\temptablewidth}{@{\extracolsep{\fill}}lllllll}
 $\tau$     &\multicolumn{2}{c}{$\alpha=0.3$}&\multicolumn{2}{c}{$\alpha=0.5$}
&\multicolumn{2}{c}{$\alpha=0.7$}\\
 \cline{2-3}\cline{4-5}\cline{6-7}
&$E_{\infty}(h,\tau)$&Rate1&$E_{\infty}(h,\tau)$
&Rate1&$E_{\infty}(h,\tau)$&Rate1\\\hline
$1/5$   & 1.5208e-2  & $\ast$  & 2.0989e-2  & $\ast$  &2.9859e-2 & $\ast$\\
$1/10$  & 3.7508e-3  & 2.0196  & 5.2056e-3  & 2.0115  &7.5354e-3 & 1.9864\\
$1/20$  & 9.3437e-4  & 2.0051  & 1.2874e-3  & 2.0155  &1.8780e-3 & 2.0045\\
$1/40$  & 2.3411e-4  & 1.9968  & 3.1958e-4  & 2.0102  &4.6777e-4 & 2.0053\\
$1/80$  & 5.8774e-5  & 1.9939  & 7.9583e-5  & 2.0057  &1.1666e-4 & 2.0035\\
\end{tabular*}
{\rule{\temptablewidth}{0.7pt}}
\end{center}
\end{table}

 \begin{table}[hbt!]
 \begin{center}
 \caption{Numerical convergence orders in spatial direction with $\tau=\frac1{20000}$ when $\alpha=0.5$ for Example \ref{ex2}.}
 \label{table4}
 \renewcommand{\arraystretch}{0.8}
 \def\temptablewidth{0.5\textwidth}
{\rule{\temptablewidth}{0.7pt}}
 \begin{tabular*}{\temptablewidth}{@{\extracolsep{\fill}}lll}
 $h$      & $E_\infty(h,\tau)$& Rate2\\\hline
 $\pi/2$  & 3.4342e-3    & $\ast$\\
 $\pi/4$  & 2.0348e-4    & 4.0770\\
 $\pi/8$  & 1.2502e-5    & 4.0247\\
 $\pi/16$ & 7.7904e-7    & 4.0043\\
 $\pi/32$ & 4.9832e-8    & 3.9665\\
 \end{tabular*}
 {\rule{\temptablewidth}{0.7pt}}
 \end{center}
 \end{table}


\begin{thebibliography}{99}

 \bibitem{Podlubny} I. Podlubny, Fractional Differential Equations, Academic Press, New York, 1999.

 \bibitem{Kilbas} A. Kilbas, H. Srivastava, J. Trujillo, Theory and Applications of Fractional Differential Equations, Elsevier Science and Technology, 2006.

 \bibitem{Lubich} C. Lubich, Discretized fractional calculus, SIAM J. Math. Anal. 17 (1986) 704--719.

 \bibitem{Meerschaert} M. Meerschaert, C. Tadjeran, Finite difference approximations for fractional advection-dispersion flow equations, J. Comput. Appl. Math. 172 (2004) 65--77.

 \bibitem{Yuste2}S. Yuste, Weighted average finite difference methods for fractional diffusion equations, J. Comput. Phys. 216 (2006) 264--274.

 \bibitem{Ervin} V.J. Ervin, J.P. Roop, Variational formulation for the stationary fractional advection dispersion equation, Numer. Meth. Part. Differ. Equ. 22 (2006) 558--576.

 \bibitem{SunWu} Z. Sun, X. Wu, A fully discrete difference scheme for a diffusion-wave system, Appl. Numer. Math. 56 (2006) 193--209.

 \bibitem{Lin} Y. Lin, C. Xu, Finite difference/spectral approximations for the time-fractional diffusion equation, J. Comput. Phys. 225 (2007) 1533--1552.

 \bibitem{ZhuangP} P. Zhuang, F. Liu, V. Anh, I. Turner, New solution and analytical techniques of the implicit numerical method for the anomalous subdiffusion equation, SIAM J. Numer. Anal. 46 (2) (2008) 1079--1095.

 \bibitem{LiXu} X. Li, C. Xu, A space-time spectral method for the time fractional diffusion equation, SIAM J. Numer. Anal. 47 (3) (2009) 2018--2131.

 \bibitem{Du} R. Du, W. Cao, Z. Sun, A compact difference scheme for the fractional diffusion-wave equation, Appl. Math. Model. 34 (2010) 2998--3007.

 \bibitem{Sousaa} E. Sousa, C. Li, A weighted finite difference method for the fractional diffusion equation based on the Riemann-Liouville derivative, (2011), arXiv:1109.2345 [math.NA].

 \bibitem{Tian} W. Tian, H. Zhou, W. Deng, A class of second order difference approximations for solving space fractional diffusion equations, Math. Comput. arXiv:1201.5949 [math.NA].

 \bibitem{Zhou} H. Zhou, W. Tian, W. Deng, Quasi-compact finite difference schemes for space fractional diffusion equations, J. Sci. Comput. 56 (2013) 45--66.

 \bibitem{Huang} J. Huang, Y. Tang, L. V\'{a}zquez, J. Yang, Two finite difference schemes for time fractional diffusion-wave equation, Numer. Algor. 64 (2013) 707--720.

 \bibitem{Wang} Z. Wang, S. Vong, Compact difference schemes for the modified anomalous fractional sub-diffusion equation and the fractional diffusion-wave equation, arXiv:1310.5298 [math.NA].

 \bibitem{Wang3} Z. Wang, S. Vong, A high order ADI scheme for the two-dimensional time fractional diffusion-wave equation, to appear in Int. J. Comput. Math., arXiv:1310.6627 [math.NA].

 \bibitem{Wang2} S. Vong, Z. Wang, Compact finite difference scheme for the fourth-order fractional sub-diffusion system, to appear in Adv. Appl. Math. Mech.

 \bibitem{Langlands} T. Langlands, B. Henry, The accuracy and stability of an implicit solution method for the fractional diffusion equation, J. Comput. Phys. 205 (2005) 719--736.

 \bibitem{Gao} G. Gao, Z. Sun, A compact finite difference scheme for the fractional sub-diffusion equations, J. Comput. Phys. 230 (2011) 586--595.

 \bibitem{Yuste1} S. Yuste, L. Acedo, An explicit finite difference method and a new von Neumann-type stability analysis for fractional diffusion equations, SIAM J. Numer. Anal. 42 (5) (2005) 1862--1874.

 \bibitem{Huang2} J. Huang, Y. Tang, W. Wang, J. Yang, A compact difference scheme for time fractional diffusion equation with Neumann boundary conditions, AsiaSim 2012. Springer Berlin Heidelberg, (2012) 273--284.

 \bibitem{zhaosun} X. Zhao, Z. Sun, A box-type scheme for fractional sub-diffusion equation with Neumann boundary conditions, J. Comput. Phys. 230 (2011) 6061--6074.

 \bibitem{Ren} J. Ren, Z. Sun, X. Zhao, Compact difference scheme for the fractinoal sub-diffusion equation with Neumann boundary conditions,
    J. Comput. Phys. 232 (2013) 456--467.

 \bibitem{Rensci} J. Ren, Z. Sun, Numerical algorithm with high spatial accuracy for the fractional diffusion-wave equation with Neumann boundary conditions, J. Sci. Comput. 56 (2013) 381--408.

 \bibitem{Quarteroni} A. Quarteroni, A. Valli, Numerical approximation of partial differential equations. Springer, Berlin, 1997.

 \bibitem{Zhang2} Y. Zhang, Z. Sun, X. Zhao, Compact alternating direction implicit scheme for the two-dimensional fractional diffusion-wave equation, SIAM J. Numer. Anal. 50 (2012) 1535--1555.

\end{thebibliography}
\end{document}